%
%

\documentstyle{amsppt}
\magnification=1200
\TagsOnRight
\NoBlackBoxes
\topmatter
\NoBlackBoxes
\title Examples of Domains\\ 
with Non-Compact Automorphism Groups
\endtitle
\footnote[]{{\bf Mathematics
Subject Classification:} 32A07, 32H05, 32M05 \hfill}
\footnote[]{{\bf Keywords and Phrases:} Automorphism
groups, Reinhardt domains, circular domains.\hfill}
\author Siqi Fu \ \ \ and \ \ \ A. V. Isaev \ \ \ and \ \ \ Steven G. Krantz
\endauthor 
\abstract We give an example of a bounded, pseudoconvex, circular 
domain in ${\Bbb C}^n$ for any $n\ge 3$ with 
smooth real-analytic boundary and non-compact automorphism group, which 
is not biholomorphically equivalent
 to any Reinhardt domain. We also give an analogous example in ${\Bbb C}^2$, 
where the domain is bounded, non-pseudoconvex, is not
equivalent to any Reinhardt domain,  and the boundary is smooth   
real-analytic at all points except one.   
\endabstract
\endtopmatter  
\document
\def\qed{{\hfill{\vrule height7pt width7pt
depth0pt}\par\bigskip}}

Let $D$ be a bounded or, more generally, a hyperbolic
domain in ${\Bbb C}^n$. Denote by $\text{Aut}(D)$ the
group of biholomorphic self-mappings of $D$. The group
$\text{Aut}(D)$, with the topology given by uniform 
convergence on
compact subsets of $D$, is in fact a Lie group \cite{Kob}.

A domain $D$ is called Reinhardt if the standard action of 
the $n$-dimensional
torus ${\Bbb T}^n$ on ${\Bbb C}^n$,  
$$
z_j\mapsto e^{{i\phi}_j}z_j,\qquad {\phi}_j\in {\Bbb R},\quad
j=1,\dots,n,
$$
leaves $D$ invariant. For certain classes of domains with 
non-compact
 automorphism groups, Reinhardt domains serve as standard models up 
to biholomorphic equivalence (see e.g. \cite{R}, \cite{W}, 
\cite{BP}, \cite{GK1}, \cite{Kod}). 

It is an intriguing question whether {\it any} domain in 
${\Bbb C}^n$ with non-compact
 automorphism group and satisfying some natural geometric 
conditions is biholomorphically 
equivalent to a Reinhardt domain. The history of the study 
of domains with non-compact
 automorphism groups shows that there were expectations that 
the answer to this question
 would be positive (see \cite{Kra}). In this note we give 
examples that show that the answer 
is in fact negative.

While the domain that we shall consider in Theorem 1 
below has already been noted in the literature \cite{BP}, 
it has never been proved that 
this domain is not biholomorphically equivalent to a 
Reinhardt domain. 
Note that this domain is circular, i.e. 
it is invariant under the special rotations
$$
z_j\mapsto e^{{i\phi}}z_j,\qquad {\phi}\in {\Bbb R},\quad
j=1,\dots,n.
$$
Our first result is the following

\proclaim{Theorem 1} There exists a bounded, pseudoconvex, 
circular domain $\Omega\subset{\Bbb C}^3$
 with smooth real-analytic boundary and non-compact 
automorphism group, which is not 
biholomorphically equivalent to any Reinhardt domain.
\endproclaim
\demo{Proof} Consider the domain
$$
\Omega=\{|z_1|^2+|z_2|^4+|z_3|^4+(\overline{z_2}z_3+
\overline{z_3}z_2)^2<1\}.
$$
The domain $\Omega$ is invariant under the action of 
the two-dimensional torus ${\Bbb T}^2$
$$
\aligned
&z_1\mapsto e^{{i\phi}_1}z_1,\qquad {\phi}_1\in 
{\Bbb R},\\
&z_j\mapsto e^{{i\phi}_2}z_j,\qquad {\phi}_2\in 
{\Bbb R}, \quad j=2,3,
\endaligned 
$$
and therefore is circular. It is also a pseudoconvex, 
bounded domain with smooth real-analytic boundary. 
The automorphism group $\text{Aut}(\Omega)$  is 
non-compact since it contains the following subgroup
$$
\aligned
&z_1\mapsto\frac{z_1-a}{1-\overline{a}z_1},\\
&z_2\mapsto\frac{(1-|a|^2)^{\frac{1}{4}}z_2}
{(1-\overline{a}z_1)^{\frac{1}{2}}},\\
&z_3\mapsto\frac{(1-|a|^2)^{\frac{1}{4}}z_3}
{(1-\overline{a}z_1)^{\frac{1}{2}}},
\endaligned \tag{1}
$$
for a complex parameter $a$ with $|a|<1$.

We are now going to explicitly determine 
$\text{Aut}(\Omega)$. Let $F=(f_1,f_2,f_3)$ be an
 automorphism of $\Omega$. Then, since $\Omega$ is 
bounded, pseudoconvex and has real-analytic
 boundary, $F$ extends smoothly to 
$\overline{\Omega}$ \cite{BL}.
 Therefore, $F$ must preserve
 the rank of the Levi form ${\Cal L}_{\partial\Omega}(q)$ 
of $\partial\Omega$ at every $q\in\partial\Omega$.
 The only points where ${\Cal L}_{\partial\Omega}\equiv 0$  
are those of the form $(e^{i\alpha},0,0)$,
 $\alpha\in{\Bbb R}$.
These points must be preserved by $F$.
 This observation implies that $f_j(e^{i\alpha},0,0)=0$ for 
all $\alpha\in{\Bbb R}$, $j=2,3$.
 Restricting $f_2$, $f_3$ to the unit disc 
$\Omega\cap\{z_2=z_3=0\}$, we see that $f_j(z_1,0,0)=0$ for
 all $|z_1|\le 1$, $j=2,3$. Therefore, $F(0)=(b,0,0)$ for some $|b|<1$.
Taking the composition of $F$
 and the automorphism $G$ of the form (1) with $a=b$, 
we find that the mapping $G\circ F$ preserves the
 origin. Since $\Omega$ is circular, it follows from a theorem of 
H. Cartan \cite{C} that $G\circ F$ must be linear.
 Therefore, any automorphism of $\Omega$ is the composition of a linear automorphism and an automorphism of the form (1).

The above argument also shows that any linear 
automorphism of $\Omega$ can be written as
$$
\aligned
&z_1\mapsto e^{i\phi_1}z_1,\\
&z_2\mapsto az_2+bz_3,\\
&z_3\mapsto cz_3+dz_3,
\endaligned
$$
where $\phi_1\in{\Bbb R}$, $a,b,c,d\in {\Bbb C}$, 
and the transformation in the variables $(z_2, z_3)$ is
 an automorphism of the section $\Omega\cap \{z_1=0\}$. 
Further, since the only points of $\partial\Omega$
 where $\text{rank}\,{\Cal L}_{\partial\Omega}=1$  are 
those of the form $(z_1,w,\pm w)$ with $w\ne 0$ and since automorphisms of
 $\Omega$ preserve such points, it follows that any 
linear automorphism of $\Omega$ is in fact given by
$$
\aligned
&z_1\mapsto e^{i\phi_1}z_1,\\
&z_2\mapsto e^{i\phi_2}z_{\sigma(2)},\\
&z_3\mapsto \pm e^{i\phi_2}z_{\sigma(3)},
\endaligned
$$
where $\phi_1, \phi_2\in{\Bbb R}$, and $\sigma$ is a 
permutation of the set $\{2,3\}$.

The preceding description of $\text{Aut}(\Omega)$ 
implies that $\text{dim}\,\text{Aut}(\Omega)=4$.
That is to say, each of the four connected components of $\text{Aut}(\Omega)$  is parametrized by the point $a$ from the unit disc and by the 
rotation parameters $\phi_1, \phi_2$.

Suppose now that $\Omega$ is biholomorphically 
equivalent to a Reinhardt domain $D\subset{\Bbb C}^3$.
 Since $\Omega$ is bounded, it follows that $D$ is hyperbolic. 
It follows from \cite{Kru} that any hyperbolic Reinhardt
 domain $G\subset{\Bbb C}^n$ can be biholomorphically 
mapped onto its normilized form $\tilde G$ for which the identity
 component $\text{Aut}_0(\tilde G)$ of $\text{Aut}(\tilde G)$ is 
described as follows. There exist integers $0\le s\le
t\le p\le n$ and $n_i\ge 1$, $i=1,\dots,p$, with $\sum_{i=1}^p
n_i=n$, and real numbers ${\alpha}_i^k$, $i=1,\dots,s$,
$k=t+1,\dots,p$, and  ${\beta}_j^k$, $j=s+1,\dots,t$,
$k=t+1,\dots,p$  such that if we set
$z^i=\left(z_{n_1+\dots+n_{i-1}+1},\dots,
z_{n_1+\dots+n_i}\right)$,
$i=1,\dots,p$, then $\text{Aut}_0(\tilde G)$ is given by 
the mappings
$$
\aligned
&z^i\mapsto\frac{A^iz^i+b^i}{c^iz^i+d^i},
\quad i=1,\dots,s,\\
&z^j\mapsto B^jz^j+e^j,\quad j=s+1,\dots,t,\\
&z^k\mapsto
C^k\frac{\prod_{j=s+1}^t\exp\left(-\beta_j^k
\left(2\overline{e^j}^TB^jz^j+
|e^j|^2\right)\right)z^k}{\prod_{i=1}^s
(c^iz^i+d^i)^{2\alpha_i^k}},\quad
k=t+1,\dots,p,
\endaligned\tag{2}
$$
where
$$
\aligned
&\pmatrix
A^i&b^i\\
c^i&d^i
\endpmatrix\in SU(n_i,1),\quad i=1,\dots,s,\\
&B^j\in U(n_j),\quad e^j\in{\Bbb C}^{n_j},
\quad j=s+1,\dots,t,\\
&C^k\in U(n_k),\quad k=t+1,\dots,p.
\endaligned
$$
The normalized form $\tilde G$ is written as
$$
\aligned
G=\Biggl\{\left|z^1\right|&<1,\dots,\left|z^s\right|<1,\\
&\Biggl(\frac{z^{t+1}}{\prod_{i=1}^s
\left(1-\left|z^i\right|^2\right)^{{\alpha}_i^{t+1}}\prod_{j=s+1}^t
\exp\left(-{\beta}_j^{t+1}\left|z^j\right|^2\right)}\ ,\ \dots \ , \\
&\frac{z^{p}}{\prod_{i=1}^s
\left(1-\left|z^i\right|^2\right)^{{\alpha}_i^p}\prod_{j=s+1}^t
\exp\left(-{\beta}_j^p\left|z^j\right|^2\right)}\Biggr)\in\tilde
G_1\Biggr\},
\endaligned\tag{3}
$$
where $\tilde G_1:=\tilde G\bigcap\left\{z^i=0,\,
i=1,\dots,t\right\}$ is a hyperbolic Reinhardt domain in ${\Bbb
C}^{n_{t+1}}\times\dots\times{\Bbb C}^{n_p}$.

It is now easy to see that, for any hyperbolic 
Reinhardt domain $D\subset{\Bbb C}^3$ written in a 
normilized form $\tilde D$,
 $\text{Aut}_0(\tilde D)$ given by formulas (2) cannot 
have dimension equal to 4.

This completes the proof.\qed
\enddemo

{\bf Remark.} The theorem can be easily 
extended to ${\Bbb C}^n$ for any $n\ge 3$ 
(just replace $|z_1|^2$ in
 the defining function of $\Omega$ by 
$\sum_{j=1}^{n-2}|z_j|^2$, $z_2$ by $z_{n-1}$, $z_3$ by $z_n$). 

There is considerable evidence that, in complex dimension two, an example
such as that constructed in Theorem 1
does not exist.  Certainly the example provided above 
depends on the decoupling,
in the domain $\Omega$, of the variables $z_2, z_3$ from the
variable $z_1$.  Such decoupling is not possible when the dimension
is only two.

The work of 
Bedford and Pinchuk (see \cite{BP} and references therein)
suggests that the only 
smoothly bounded domains in ${\Bbb C}^2$ with non-compact automorphism groups
are (up to biholomorphic equivalence) the complex ellipsoids
$$
\Omega_{\alpha} = \{(z_1,z_2) \in 
{\Bbb C}^2: |z_1|^2 + |z_2|^{2\alpha} < 1\},
$$
where $\alpha$ is a positive integer. It is also a plausible 
conjecture that any bounded domain in ${\Bbb C}^2$ 
with non-compact automorphism group and a boundary of finite smoothness $C^k$ for $k\ge 1$, 
is biholomorphically equivalent to some $\Omega_{\alpha}$, where $\alpha\ge 1$ and is not necessarily an integer.   
Of course all the domains $\Omega_{\alpha}$ are pseudoconvex and Reinhard.

However, as the following theorem shows, if we allow the boundary to be non-smooth at just one point, then
the domain may be non-pseudoconvex and be non-equivalent to any Reinhardt domain.

\proclaim{Theorem 2} There exists a bounded, non-pseudoconvex domain $\Omega\subset{\Bbb C}^2$ 
with non-compact automorphism group such that $\partial\Omega$ 
is smooth real-analytic everywhere except one point (this
exceptional point is an orbit accumulation point for the automorphism group action), 
and such that $\Omega$ is not biholomorphically equivalent to any Reinhardt domain.
\endproclaim

For the proof of Theorem 2, we first need the following lemma.

\proclaim{Lemma A} If $\Omega\subset{\Bbb C}^2$ is a bounded, non-pseudoconvex, simply-connected domain such that the identity component $\text{Aut}_0(\Omega)$ of the automorphism group $\text{Aut}(\Omega)$ is non-compact, 
then $\Omega$ is not biholomorphically equivalent to any Reinhardt domain.
\endproclaim

\demo{Proof of Lemma A} Suppose that $\Omega$ is biholomorphically equivalent to a Reinhardt domain $D$. 
Since $\Omega$ is bounded, it follows that $D$ is hyperbolic. Also, since $\text{Aut}_0(\Omega)$ is 
non-compact, then so is $\text{Aut}_0(D)$. We are now going to show that any such domain $D$ is either 
pseudoconvex, or not simply-connected, or cannot be biholomorphically equivalent to a bounded domain. 
This result clearly implies the lemma.

We can now assume that the domain $D$ is written in its normalized form $\tilde D$ as in (3), and $\text{Aut}_0(\tilde D)$ is given by formulas (2). Then, since $\text{Aut}_0(\tilde D)$ is non-compact, it must be that $t>0$. Next, 
if $p=t$, then $\tilde D$ is either non-hyperbolic (for $s<t$), or (for $s=t$) is the unit ball or the unit polydisc and therefore is pseudoconvex. Thus we can assume that $t=1$, $p=2$, $n_1=n_2=1$.

Let $\tilde D_1\subset{\Bbb C}$ be the hyperbolic Reinhardt domain analogous to $\tilde G_1$ that was defined above (see (3)).
Clearly, there are the following possibilities for $\tilde D_1$:
\medskip 

\noindent {\bf (i)} $\tilde D_1=\{0<|z_2|<R\}$, $0<R<\infty$;

\noindent {\bf (ii)} $\tilde D_1=\{r<|z_2|<R\}$, $0<r<R\le \infty$;

\noindent {\bf (iii)} $\tilde D_1=\{|z_2|<R\}$, $0<R<\infty$.
\medskip 

For the cases {\bf (i)}, {\bf (ii)}, $\tilde D$ is always not simply-connected, and therefore we will concentrate on the case {\bf (iii)}.
 If $s=0$, then $\tilde D$ is not hyperbolic since it contains the complex line $\{z_2=0\}$.  
Thus we can assume that $s=1$. Next observe that, for $\alpha_1^2\ge 0$, 
the domain $\tilde D$ is always pseudoconvex.  Thus we may take $\alpha_1^2<0$. Then the domain $\tilde D$ has the form
$$
\tilde D=\left\{|z_1|<1,\,|z_2|<\frac{R}{(1-|z_1|^2)^{\gamma}}\right\},\qquad \gamma>0.
$$
We will now show that the above domain $\tilde D$ cannot be biholomorphically equivalent 
to a bounded domain. More precisely, we will show that any bounded holomorphic 
function on $\tilde D$ is independent of $z_2$.

Let $f(z_1,z_2)$ be holomorphic on $\tilde D$ and $|f|<M$ for some $M>0$. For every $\rho$ such that $|\rho|\le\frac{R}{2}$, the disc $\Delta_{\rho}=\{|z_1|<1,\,z_2=\rho\}$ is contained in $\tilde D$. We will show that ${\partial f}/\partial z_2 \equiv 0$
  on every such $\Delta_{\rho}$, which implies that ${\partial f}/\partial z_2 \equiv 0$ everywhere in $\tilde D$.

Fix a point $(\mu,\rho)\in \Delta_{\rho}$ and restrict $f$ to the disc $\Delta'_{\mu}=\{z_1=\mu,\,|z_2|<R_{\mu}\}$, where $R_{\mu}={R}/{2(1-|\mu|^2)^{\gamma}}$ . Clearly, $(\mu,\rho)\in\Delta'_{\mu}$ and $\overline{\Delta'_{\mu}}\subset \tilde D$. By the 
Cauchy Integral Formula
$$
f(\mu,z_2)=\frac{1}{2\pi i}\int_{\partial \Delta'_{\mu}}\frac{f(\mu,\zeta)}{\zeta-z_2}\,d\zeta,
$$
for $|z_2|<R_{\mu}$, and therefore
$$
\frac{\partial f}{\partial z_2}(\mu,\rho)=\frac{1}{2\pi i}\int_{\partial \Delta'_{\mu}}\frac{f(\mu,\zeta)}{(\zeta-\rho)^2}\,d\zeta.
$$
Hence
$$
\left|\frac{\partial f}{\partial z_2}(\mu,\rho)\right|\le\frac{MR_{\mu}}{(R_{\mu}-|\rho|)^2}.
$$
Letting $|\mu|\rightarrow 1$ and taking into account that $R_{\mu}\rightarrow\infty$, we see that
$\left|{\partial f}/\partial z_2 (\mu,\rho)\right|\rightarrow 0$ as  $|\mu|\rightarrow 1$. Therefore, ${\partial f}/\partial z_2 \equiv 0$ on $\Delta_{\rho}$.

The lemma is proved.\qed
\enddemo

\demo{Proof of Theorem 2} We will now present a domain that satisfies the conditions of the lemma. Set
$$
\Omega=\left\{|z_1|^2+|z_2|^4+8|z_1-1|^2\left(\frac{z_2^2}{z_1-1}-\frac{3}{2}\frac{|z_2|^2}{|z_1-1|}+\frac{\overline{z_2}^2}{\overline{z_1}-1}\right)^2<1\right\}.
$$

The domain $\Omega$ is plainly bounded since the third term on the left is non-negative. Next, the identity component $\text{Aut}_0(\Omega)$ of its automorphism group is non-compact since it contains the subgroup
$$
\aligned
&z_1\mapsto\frac{z_1-a}{1-az_1},\\
&z_2\mapsto\frac{(1-a^2)^{\frac{1}{4}}z_2}
{(1-az_1)^{\frac{1}{2}}},
\endaligned
$$
where $a\in(-1,1)$. 

Further, $\Omega$ is simply-connected, since the family of mappings $F_{\tau}(z_1,z_2)=(z_1,\tau z_2)$, $0\le\tau\le 1$, 
retracts $\Omega$ inside itself, as $\tau\rightarrow 0$, to the unit disc $\{|z_1|<1,\, z_2=0\}$ (which is simply-connected).

To show that $\Omega$ is not pseudoconvex, consider its unbounded realization. Namely, under the mapping
$$
\aligned
&z_1\mapsto \frac{z_1+1}{z_1-1},\\
&z_2\mapsto \frac{\sqrt{2}z_2}{\sqrt{z_1-1}},
\endaligned\tag{4}
$$
the domain $\Omega$ is transformed into the domain
$$
\Omega'=\left\{\text{Re}\,z_1+\frac{1}{4}|z_2|^4+2\left(z_2^2-\frac{3}{2}|z_2|^2+\overline{z_2}^2\right)^2<0\right\}.
$$
It is easy to see that at the boundary point $(-\frac{3}{4},1)\in\partial\Omega'$ the Levi form of $\partial\Omega'$ is equal to $-|z_2|^2$, and thus is negative-definite. Therefore, $\Omega$ is non-pseudoconvex.

Hence, by Lemma A, $\Omega$ is not biholomorphically equivalent to any Reinhardt domain. 

Next, if $\phi$ denotes the defining function of $\Omega$, the following holds at every boundary point 
of $\Omega$ except $(1,0)$:
$$
\frac{\partial \phi}{\partial z_1}=\frac{1}{z_1-1}\left(-\frac{z_2}{2}\frac{\partial \phi}{\partial z_2}+1-\overline{z_1}\right),
$$
and therefore $\text{grad}\,\phi$ does not vanish at every such point. Hence, $\partial\Omega$ is smooth real-analytic everywhere except at $(1,0)$.

The theorem is proved.\qed
\enddemo

{\bf Remarks.}   \break
\smallskip

\noindent {\bf 1.}  The hypothesis of 
simple connectivity in Lemma A is automatically satisfied if, for example, the boundary of the domain is locally variety-free and smooth near some orbit accumulation point for the automorphism group of the domain (see e.g. \cite{GK2}). For a smoothly boun
ded domain it would follow from a conjecture of Greene/Krantz \cite{GK3}.

\noindent {\bf 2.}  Tedious calculations show that the 
boundary of the domain $\Omega$ in Theorem 2 is 
quite pathological near the exceptional point $(1,0)$.  
It is not Lipschitz-smooth of any positive degree.  It would be
interesting to know whether there is an example with Lipschitz-1
boundary at the bad point. 

In fact, many more examples similar to that in Theorem 2 can be
constructed in the following way.  Let 
$$
\Omega' = \{(z_1,z_2) \in {\Bbb C}^2: 
\hbox{Re}\, z_1 + P(z_2) < 0\},\tag{5}
$$ 
where $P=|z_2|^{2m}+Q(z_2)$ is a homogeneous 
non-plurisubharmonic polynomial, $m$ is a positive 
integer, and $Q(z_2)$ is positive away from the origin. 
Then, by a mapping analogous to (4), $\Omega'$ can be transformed
into a bounded domain $\Omega$.  The domain $\Omega$ is simply-connected,
non-pseudoconvex, $\text{Aut}_0(\Omega)$ is non-compact, and
$\partial \Omega$ is smooth real-analytic everywhere except at the
point $(1,0)$.  For all such examples, $\partial \Omega$
is not Lipschitz-smooth of any positive degree at $(1,0)$.

It is also worth noting that, in the example contained in Theorem 2, the
point $(-1,0)$ is also an orbit accumulation point, but $\partial \Omega$
is smooth real-analytic at this point. 

\noindent {\bf 3.} It is conceivable that
the domain $\Omega$ as in Theorem 2 has an alternative, smoothly bounded
realization, but it looks plausible that if in formula (5) we allow $P(z_2)$ to be an arbitrary homogeneous polynomial positive away from the origin with no harmonic term, then domain (5) does not have a bounded realization with $C^1$-smooth boundary, unl
ess $P(z_2)=c|z_2|^{2m}$, where $c>0$ and $m$ is a positive integer. 
\bigskip 

This work was completed while the second 
author was an Alexander von Humboldt Fellow at the 
University of Wuppertal.  Research at MSRI by the third author was
supported in part by NSF Grant DMS-9022140.

\bigskip
\pagebreak

{\obeylines
Department of Mathematics 
University of California, Irvine, CA 92717
USA 
E-mail address: sfu\@math.uci.edu
\bigskip

Centre for Mathematics and Its Applications 
The Australian National University 
Canberra, ACT 0200
AUSTRALIA 
E-mail address: Alexander.Isaev\@anu.edu.au
\smallskip
and
\smallskip
Bergische Universit\"at
Gesamthochschule Wuppertal
Mathematik (FB 07)
Gaussstrasse 20
42097 Wuppertal
GERMANY
E-mail address: Alexander.Isaev\@math.uni-wuppertal.de
\bigskip

Department of Mathematics
Washington University, St.Louis, MO 63130
USA 
E-mail address: sk\@math.wustl.edu
\smallskip
and
\smallskip
MSRI
1000 Centennial Drive
Berkeley, California 94720
USA
E-mail address: krantz\@msri.org
\smallskip

\enddocument